\documentclass[aop,seceqn,MSNbibl,dvips]{arximspdf}
\usepackage{mathbh}

\doi{10.1214/09-AOP488}
\volume{38}
\issue{2}
\pubyear{2010}
\firstpage{479}
\lastpage{497}

\makeatletter

\DeclareMathAlphabet\mathcaligr{OMS}{cmsy}{m}{n}

\newtheorem{theorem}{Theorem}
\newtheorem{proposition}{Proposition}[section]

\newtheorem{lemma}[proposition]{Lemma}
\newtheorem{corollary}[proposition]{Corollary}

\makeatother

\begin{document}
\begin{frontmatter}

\title{A Trotter-type approach to infinite rate mutually catalytic branching\protect\thanksref{T1}}
\runtitle{A trotter type approach to IMUB}

\begin{aug}
\author[A]{\fnms{Achim} \snm{Klenke}\ead[label=e1]{math@aklenke.de}\corref{}}
\and
\author[A]{\fnms{Mario} \snm{Oeler}\ead[label=e2]{Mario.Oeler@web.de}}
\thankstext{T1}{Supported in part by the German Israeli Foundation,
Grant number G-807-227.6/2003.}
\runauthor{A.~Klenke and M.~Oeler}
\affiliation{Johannes Gutenberg-Universit{\"a}t Mainz}
\address[A]{Institut f{\"u}r
Mathematik\\Johannes Gutenberg-Universit{\"a}t Mainz\\Staudingerweg 9\\55009 Mainz\\Germany\\
\printead{e1}\\
\phantom{E-mail:\ }\printead*{e2}} 
\end{aug}

\received{\smonth{4} \syear{2009}}

%
\begin{abstract}
Dawson and Perkins [\textit{Ann. Probab.} \textbf{26} (1988)
1088--1138] constructed a stochastic model
of an interacting two-type population indexed by a countable site space
which locally undergoes a mutually
catalytic branching mechanism. In Klenke and Mytnik [Preprint (2008),
 arXiv:\href{http://arxiv.org/abs/0901.0623}{0901.0623}], it is shown
that as the branching rate approaches infinity, the process converges
to a process that is called the
\emph{infinite rate mutually catalytic branching process} (IMUB). It is
most conveniently characterized
as the solution of a certain martingale problem. While in the latter
reference, a noise equation approach
is used in order to construct a solution to this martingale problem,
the aim of this paper is to provide a Trotter-type construction.

The construction presented here will be used in a forthcoming paper,
Klenke and Mytnik [Preprint (2009)],
to investigate the long-time behavior of IMUB (coexistence versus
segregation of types).

This paper is partly based on the Ph.D. thesis of the second author
(2008), where the Trotter approach was first introduced.
\end{abstract}

%


\end{frontmatter}

\section{Introduction and main results}
\label{S1}
\subsection{Background and motivation}
\label{S1.1}

In~\cite{DawsonPerkins1998}, Dawson and Perkins studied a stochastic
model of mutually catalytic (continuous-state) branching. Two
populations live on a countable site space $S$ and the amount of
population of type $i=1,2$ at time $t$ at site $k\in S$ is denoted by
$Y_{i,t}(k)\in[0,\infty)$. The populations migrate according to a
deterministic heatflow-like dynamics that is characterized by the
(symmetric) $q$-matrix ${\mathcaligr{A}}$ of a Markov chain on $S$.
Locally, the
populations undergo critical continuous-state branching with a rate
that is proportional to the size of the other type at the same place.
Formally, this model can be described by a system of stochastic
differential equations:
\begin{eqnarray}
\label{E1.1}
\ Y_{i,t}(k)&=&Y_{i,0}(k)+\int_0^t \sum_{l\in S} {\mathcaligr{A}}(k,l)
Y_{i,s}(l)
\, ds\nonumber
\\[-8pt]\\[-8pt]
&&{}+\int_0^t (\gamma Y_{1,s}(k)Y_{s,2}(k))^{1/2}
\,dW_{i,s}(k),\qquad   t\geq0, k\in S, i=1,2.\nonumber
\end{eqnarray}
Here, $(W_i(k),  k\in S, i=1,2)$ is an independent family of
one-dimensional Brownian motions and $Y_0$ is chosen from a suitable
subspace of $([0,\infty)^2)^S$. The parameter $\gamma\geq0$ can be
thought of as being the branching rate for this model.
Dawson and Perkins showed that there is a unique weak solution of \mbox
{\rm(\ref{E1.1})} and studied the long-time behavior of this
model. They also constructed the analogous model in the continuous
setting on ${\mathbb{R}}$ instead of $S$.

For the model with $S={\mathbb{Z}}$ and ${\mathcaligr{A}}$ the
$q$-matrix of symmetric nearest
neighbor random walk, the model tends to a state with spatially
segregated types. In an approach to describing the cluster growth
quantitatively, a space and time rescaling argument suggests that it is
useful to first study the limit as $\gamma\to\infty$. Studying this
limit requires a formal description of the limit process $X$,
construction of the limit process and the establishing of convergence
of $Y$ as $\gamma\to\infty$.

This program is carried out for a process where $S$ is a singleton in
\cite{KM1} and for a countable site space $S$ in \cite{KM2}.
Furthermore, in \cite{KM3}, the long-time behavior is studied which
shows a dichotomy between coexistence and segregation of types,
depending on the potential properties of the matrix ${\mathcaligr{A}}$.

In \cite{KM2}, the process $X$ is characterized both via a martingale
problem and as the solution of a system of stochastic differential
equations of jump type. While the construction of $X$ was performed via
the construction of approximate solutions of the stochastic
differential equations, here, the aim is to present a different
approach via a Trotter approximation scheme.

The main idea is described via the following heuristics. Denote by
$a_t$ the matrix of time $t$ transition probabilities of the
continuous-time Markov chain with $q$-matrix ${\mathcaligr{A}}$.
Furthermore, let
$Q_t(y,dy')$ denote the transition kernel for equation~\mbox{\rm(\ref
{E1.1})} with
${\mathcaligr{A}}=0$. It is not hard to see that $Q_t$ converges, as
$t\to\infty$,
to some kernel~${\mathbf{Q}}$. In fact, if ${\mathcaligr{A}}=0$, then
all colonies evolve
independently and each colony is a time-transformed planar Brownian
motion in $(0,\infty)^2$, stopped when it hits the boundary. Hence,
${\mathbf{Q}}
$ is the product of the harmonic measures of planar Brownian motions in
the upper-right quadrant. Now, let ${\varepsilon }>0$, define
$X^{\varepsilon }_0=Y_0$ and
inductively let
$X^{\varepsilon }_{(k+1){\varepsilon }}$ be distributed, given
$X^{\varepsilon }_{k{\varepsilon }}$, like ${\mathbf{Q}}
(a_{\varepsilon }X^{\varepsilon }_{k{\varepsilon }},dy')$. This
amounts to an interlaced dynamics where
deterministic heatflow and random infinite rate branching alternate.
The main result of this paper is that the processes $X^{\varepsilon }$
in fact
converge, as ${\varepsilon }\to0$, to the infinite rate mutually catalytic
branching process $X$ constructed in \cite{KM2}. In Sections~\ref{S1.2}
and \ref{S1.3}, we provide a formal description of this~$X$.

The idea of using a Trotter-type approach for the construction of the
infinite rate mutually catalytic branching process is taken from the
Ph.D. thesis \cite{Oeler2008} and parts of the strategy of proof are
based on that thesis.

While the noise equation approach of \cite{KM2} relies on a duality of
the processes in order to show convergence of a sequence of
approximating processes, the Trotter approach works without this
duality. This greater flexibility is exploited in~\cite{KM3} for the
construction of a process $X^K$ with state space $([0,K]^2\setminus
(0,K)^2)^S$ that approaches $X$ and whose coordinate processes are
driven by orthogonal $L^2$-martingales. For this process $X^K$ that is
used in order to study the long-time behavior of $X$, we do not have a
duality and, thus, the noise equation approach does not seem to be feasible.

Furthermore, we hope that the Trotter-type approach could serve as a
key tool for the construction of infinite rate symbiotic branching
processes. Symbiotic branching processes with index $\varrho\in[-1,1]$
are solutions of \mbox{\rm(\ref{E1.1})}, but with $W_1(k)$ and
$W_2(k)$ being
\emph
{correlated} Brownian motions with correlation $\varrho$. These were
introduced in \cite{EtheridgeFleischmann2004}. Clearly, $\varrho=0$ is
the branching case considered here, $\varrho=-1$ is the case of
interacting Wright--Fisher diffusions and $\varrho=1$ is the parabolic
Anderson model. The voter model can be considered as the infinite rate
interacting Wright--Fisher diffusion model and can be obtained rather
simply from this model via the Trotter approach. The other cases of
$\varrho$ are open. For symbiotic branching processes, there is a
moment dual, but it is of limited use in many cases. Hence, the
Trotter-type approach might also prove useful here to construct
infinite rate versions of these processes.

\subsection{The infinite rate branching process}
\label{S1.2}
We start with a definition of the state spaces of our processes.
Define
$E:=[0,\infty)^2\setminus(0,\infty)^2$.
Let $S$ be a countable set.
For $u,v\in[0,\infty)^S$, define
\[
\langle u,v\rangle=\sum_{k\in S}u(k)v(k) \in[0,\infty].
\]
Similarly, for $x\in([0,\infty)^2)^S$ and $\zeta\in[0,\infty)^S$, define
\[
\langle x,\zeta\rangle=\sum_{k\in S}x(k)\zeta(k) \in[0,\infty]^2.
\]

We can weaken the requirement that ${\mathcaligr{A}}$ be a $q$-matrix:
let ${\mathcaligr{A}}=({\mathcaligr{A}}(k,l))_{k,l\in S}$ be a matrix
indexed by the countable
set $S$ satisfying
\begin{equation}
\label{E1.2}
{\mathcaligr{A}}(k,l)\geq0\qquad\mbox{for }k\not= l
\end{equation}
and
\begin{equation}
\label{E1.3}
\|{\mathcaligr{A}}\|:=\sup_{k \in S}\sum_{l\in S} |{\mathcaligr{A}}(k,l)|
+ |{\mathcaligr{A}}(l,k)|
<\infty.
\end{equation}

By Lemma~IX.1.6 of \cite{Liggett1985}, there exists a $\beta\in
(0,\infty
)^S$ and an $M\geq1$ such that
$\sum_{k\in S}\beta(k)< \infty$,
and
\begin{equation}
\label{E1.4}
\sum_{l\in S}\beta(l)\bigl(|{\mathcaligr{A}}(k,l)|+|{\mathcaligr
{A}}(l,k)|\bigr)\leq M\beta(k)\qquad\mbox{for all }
k\in S.
\end{equation}
We fix this $\beta$ for the rest of the paper.

Define the spaces
\[
\begin{aligned}
\mathbb{L}^{\beta}&=&\{ u\in[0,\infty)^S\dvtx \langle u,\beta
\rangle<\infty
\},\\
\mathbb{L}^{\beta,2}&=&\{ x\in([0,\infty)^2)^S\dvtx
\langle x,\beta
\rangle\in
[0,\infty)^2\},\\
\mathbb{L}^{f,2}&=&\{ y\in([0,\infty)^2)^S\dvtx y (k)\neq0\mbox{
for only finitely many }k\in S\},
\end{aligned}
\]
as well as
\[
\mathbb{L}^{\beta,E}=\mathbb{L}^{\beta,2}\cap E^S\quad\mbox
{and}\quad\mathbb{L}^{f,E}=\mathbb{L}^{f,2}\cap E^S.
\]

Finally, define the spaces
\begin{eqnarray}
\label{E1.5}
\mathbb{L}_\infty^{\beta}&=&\{f\in[0,\infty)^S\dvtx \langle
f,g\rangle<\infty\mbox{ for all }
g\in\mathbb{L}^{\beta}\}\nonumber
\\[-8pt]\\[-8pt]
&=&\Bigl\{f\in\mathbb{L}^{\beta}\dvtx \sup_{k\in S}f(k)/\beta(k)<\infty
\Bigr\}\nonumber
\end{eqnarray}
and
\[
\mathbb{L}_\infty^{\beta,E}=\{ \eta=(\eta_1,\eta_2)\in E^S\dvtx
\eta_1, \eta_2\in
\mathbb{L}_\infty^{\beta}
\}.
\]

Let ${\mathcaligr{A}}f(k)=\sum_{l\in S}{\mathcaligr{A}}(k,l)f(l)$ if the
sum is well defined.
Let ${\mathcaligr{A}}^n$ denote the $n$th matrix power of ${\mathcaligr
{A}}$ [note that this is
well defined and finite by \mbox{\rm(\ref{E1.3})}] and define
\[
a_t(k,l):=e^{t{\mathcaligr{A}}}(k,l):=\sum_{n=0}^\infty\frac{t^n
{\mathcaligr{A}}^n(k,l)}{n!}.
\]
Let $\mathcaligr{S}$ denote the (not necessarily Markov) semigroup
generated by
${\mathcaligr{A}}$, that is,
\[
\mathcaligr{S}_t f(k)= \sum_{l\in S} a_t(k,l)f(l)\qquad\mbox{for }t\geq0.
\]
We will also use the notation ${\mathcaligr{A}}f$, $\mathcaligr{S}_tf$ and
so on for
$[0,\infty
)^2$-valued functions $f$ with the obvious coordinate-wise meaning.

For $u\in{\mathbb{R}}^S$, define
\begin{equation}
\label{E1.6}
\| u\|_{\beta} =\sum_{k\in S} |u(k)|\beta(k).
\end{equation}
Note that for $f\in\mathbb{L}^\beta$, the expressions ${\mathcaligr
{A}}f$ and $\mathcaligr{S}_t f$
are well defined and that [recall $M$ from \mbox{\rm(\ref{E1.4})}]
\begin{equation}
\label{E1.7}
\Vert{\mathcaligr{A}}f\Vert_\beta\leq M\Vert f\Vert_\beta
\quad\mbox{and}\quad
\Vert\mathcaligr{S}_t f\Vert_\beta\leq e^{Mt}\Vert f\Vert_\beta.
\end{equation}
That is, the spaces $\mathbb{L}^{\beta}$ and $\mathbb{L}^{\beta,2}$
are preserved under the dynamics
of $(\mathcaligr{S}_t)$.

Let $D([0,\infty);\mathbb{L}^{\beta,E})$ be the Skorohod space of
c{\`a}dl{\`a}g
$\mathbb{L}^{\beta,E}
$-valued functions.

We will employ a martingale problem in order to characterize the
infinite rate mutually catalytic branching process $X\in D([0,\infty
);\mathbb{L}^{\beta,E})$.
In order to formulate this martingale problem for $X$ conveniently, for
$x=(x_1,x_2)\in{\mathbb{R}}^2$ and $y=(y_1,y_2)\in{\mathbb{R}}^2$,
we introduce the
\emph{lozenge product}
\begin{equation}
\label{E1.8}
x \diamond {} y := -(x_1+x_2)(y_1+y_2) + i(x_1-x_2)(y_1-y_2)
\end{equation}
(with $i=\sqrt{-1}$) and define
\begin{equation}
\label{E1.9}
F(x, y)=\exp(x\diamond y).
\end{equation}
Note that $x\diamond y=y\diamond x$, hence $F$ is symmetric.
For $x,y\in({\mathbb{R}}^2)^S$, we write
\begin{equation}
\label{E1.10}
\langle\hspace*{-0.13em}\langle x, y\rangle\hspace
*{-0.13em}\rangle   = \sum_{k\in
S}x(k)\diamond y(k)
\end{equation}
whenever the infinite sum is well defined
and let
\begin{equation}
\label{E1.11}
H(x, y)=\exp(\langle\hspace*{-0.13em}\langle x, y\rangle\hspace
*{-0.13em}\rangle ).
\end{equation}

Note that the function $H(x,y)$ is well defined if either $x\in
({\mathbb{R}}
^2)^S$ and $y\in\mathbb{L}^{f,E}$ or $x\in\mathbb{L}^{\beta,E}$
and $y\in\mathbb{L}_\infty^{\beta,E}$.

It is shown in \cite{KM1}, Corollary 2.4, that the vector space of
finite linear combinations $\sum_{i=1}^n\alpha_iF(\bolds{\cdot
},y_i)$, $n\in
{\mathbb{N}}$,
$\alpha_i\in{\mathbb{C}}$, $y_i\in E$, is dense in the space
$C_l(E;{\mathbb{C}})$ of
bounded continuous complex-valued functions on $E$ with a limit at
infinity. Hence, the family $H(\bolds{\cdot},y)$, $y\in\mathbb
{L}^{f,E}$, is
measure-determining for probability measures on $\mathbb{L}^{\beta
,E}$ (but not on
$\mathbb{L}^{\beta,2}$).

In \cite{KM2}, the following theorem was established.\setcounter{theorem}{-1}
\begin{theorem}
\label{T0}
\textup{(a)}
For all $x\in\mathbb{L}^{\beta,E}$, there exists a unique solution
$X\in D([0,\infty
);\mathbb{L}^{\beta,E})$ of the following martingale problem:
for each $y\in\mathbb{L}^{f,E}$, the process $M^{x,y}$ defined by
\renewcommand{\theequation}{MP}
\begin{eqnarray}
\label{MP1}
M^{x,y}_t:=H(X_t,y) -H(x,y)-\int_0^t \langle\hspace
*{-0.13em}\langle{\mathcaligr{A}}X_s,y\rangle\hspace*{-0.13em}\rangle
  H(X_s,y)\, ds
\end{eqnarray}
is a martingale with $M^{x,y}_0=0$.

\textup{(b)} For any $x\in\mathbb{L}^{\beta,E}$ and $y\in\mathbb
{L}_\infty^{\beta,E}$, the process $M^{x,y}$
is well defined and is a martingale.

\textup{(c)} Denote by $P_x$ the distribution of $X$ with $X_0=x$. Then
$(P_x)_{x\in\mathbb{L}^{\beta,E}}$ is a strong Markov family.
\end{theorem}

Note that for the uniqueness, it is crucial that the single coordinates
take values in $E$. If we required only values in $[0,\infty)^2$, then
the finite rate mutually catalytic branching process $Y$ is also a
solution of the martingale problem for any $\gamma\geq0$. In
Proposition~\ref{P1.1}, we will see that our approximate process
$X^{\varepsilon }
$ is also a solution to $\mbox{\rm(\ref{MP1})}$ with the larger
state space $\mathbb{L}^{\beta,2}$.

In \cite{KM3}, Theorem 1.3, it was shown that the processes $Y$ defined
in \mbox{\rm(\ref{E1.1})} converge to $X$ as $\gamma\to\infty$ in
the Meyer--Zheng
topology. Hence, the name \emph{infinite rate mutually catalytic
branching process} for $X$ is justified.

\subsection{The main result}
\label{S1.3}
We now define the approximating process $X^{\varepsilon }$ in detail.
In order to do so, we introduce the harmonic measure $Q$ of planar Brownian
motion $B$ on $(0,\infty)^2$. That is, if $B=(B_1,B_2)$ is a Brownian
motion in ${\mathbb{R}}^2$ started at $x\in[0,\infty)^2$ and
$\tau=\inf\{t>0\dvtx  B_t\notin(0,\infty)^2\}$, then we define
\renewcommand{\theequation}{\arabic{section}.\arabic{equation}}
\setcounter{equation}{11}
\begin{equation}
\label{E1.12}
Q_x={\mathbf{P}}_x[B_\tau\in\bolds{\cdot}].
\end{equation}

Now, for fixed ${\varepsilon }>0$, consider the stochastic process
$X^{\varepsilon }$ with values
in $\mathbb{L}^{\beta,2}$ with the following dynamics:
\begin{longlist}[(ii)]
\item[(i)]
Within each time interval $[n{\varepsilon },(n+1){\varepsilon })$,
$n\in{\mathbb{N}}_0$, $X^{\varepsilon }$
is the
solution of \mbox{\rm(\ref{E1.1})} with $\gamma=0$; that is, for
$k\in S$,
\[
dX^{\varepsilon }_{i,t}(k)=({\mathcaligr{A}}X^{\varepsilon }_{i,t})(k)
\,dt\qquad\mbox{for }t\in\bigl[n{\varepsilon },(n+1){\varepsilon }\bigr).
\]
Clearly, the explicit solution is
\[
X^{\varepsilon }_{i,t}(k)=(\mathcaligr{S}_{t-n{\varepsilon }}
X^{\varepsilon }_{i,n{\varepsilon }})(k)\qquad\mbox{for }t\in
\bigl[n{\varepsilon }
,(n+1){\varepsilon }\bigr).
\]
\item[(ii)]
At time $n{\varepsilon }$, $X^{\varepsilon }$ has a discontinuity.
Independently, each
coordinate $X^{\varepsilon }_{n{\varepsilon }-}(k)=\mathcaligr
{S}_{{\varepsilon }}X^{\varepsilon }_{(n-1){\varepsilon }}(k)$ is replaced
by a random element of $E$ drawn according\vspace{1pt} to the distribution
$Q_{X^{\varepsilon }
_{n{\varepsilon }-}(k)}$.
\end{longlist}
If, for $x\in E^S$, we denote by ${\mathbf{Q}}(x,\bolds{\cdot
})=\bigotimes_{k\in
S}Q_{x(k)}$ the Markov kernel of independent displacements, then
$(X^{\varepsilon }
_{n{\varepsilon }})_{n\in{\mathbb{N}}_0}$ is a Markov chain on
$\mathbb{L}^{\beta,E}$ with transition kernel
${\mathbf{Q}}^{\varepsilon }(x,\bolds{\cdot}):={\mathbf
{Q}}(\mathcaligr{S}_{\varepsilon }x,\bolds{\cdot})$. Note that
$X^{\varepsilon }$ is a
c{\`
a}dl{\`a}g process with values in $\mathbb{L}^{\beta,2}$ (but not in
$\mathbb{L}^{\beta,E}$!) and that,
for any $y\in\mathbb{L}^{f,E}$,
\[
H(X^{\varepsilon }_t)-\int_{n{\varepsilon }}^t\langle\hspace
*{-0.13em}\langle{\mathcaligr{A}}X^{\varepsilon }
_s,y\rangle\hspace*{-0.13em}\rangle   H(X^{\varepsilon }_s,y)\, ds,
\qquad t\in\bigl[n{\varepsilon },(n+1){\varepsilon }\bigr),
\]
is a martingale. Furthermore, as we will show in Lemma~\ref{L2.1.2}, we
have\break $\int H(x',y)\times {\mathbf{Q}}(x,dx')=H(x,y)$ for all $y\in\mathbb
{L}^{f,E}$ and $x\in
\mathbb{L}^{\beta,2}
$. As an immediate consequence, we get the following proposition.

\begin{proposition}
\label{P1.1}
For all $x\in\mathbb{L}^{\beta,E}$ and $y\in\mathbb{L}^{f,E}$, and
for $X^{\varepsilon }$ defined as above
with $X_0=x$, we have that
\begin{eqnarray}
\label{E1.13}
\qquad M^{{\varepsilon },x,y}_t&:=&H(X^{\varepsilon }_t,y)-H(X^{\varepsilon
}_0,y)\nonumber
\\[-8pt]\\[-8pt]
&&{}-\int_0^t \langle
\hspace*{-0.13em}\langle{\mathcaligr{A}}X^{\varepsilon }
_s,y\rangle\hspace*{-0.13em}\rangle   H(X^{\varepsilon }_s,y)
\,ds,\qquad t\geq0,\mbox{ is a martingale. }\nonumber
\end{eqnarray}
\end{proposition}

We will show that $X^{\varepsilon }$ converges to a process that takes
values in
$\mathbb{L}^{\beta,E}$ while preserving this martingale property.

The main theorem of this paper is the following.

\begin{theorem}
\label{T1}
For any $x\in\mathbb{L}^{\beta,E}$, as ${\varepsilon }\to0$, the
processes $X^{\varepsilon }$ converge in
distribution in the Skorohod spaces $D([0,\infty);\mathbb{L}^{\beta
,2})$ to the unique
solution $X$ of the martingale problem \mbox{\rm(\ref{MP1})}.
\end{theorem}

With a small effort, this construction can be interpreted as a Trotter
product approach. Recall that (under suitable assumptions on the spaces
and cores of the operators involved), the Trotter product formula
states the following (see, e.g., \cite{EthierKurtz1986}, Corollary
6.7): if $(S_t)_{t\geq0}$, $(T_t)_{t\geq0}$ and
$(U_t)_{t\geq0}$ are strongly continuous contraction semigroups with
generators $A$, $B$ and $C=A+B$, respectively, then
\[
\lim_{{\varepsilon }\downarrow0}(T_{{\varepsilon }}S_{{\varepsilon
}})^{\lfloor t/{\varepsilon }\rfloor} =
U_t\qquad
\mbox{pointwise.}
\]
In our setting, $T_t={\mathbf{Q}}$ for all $t>0$ and $T_0=\mathrm
{id}$, hence
$(T_t)$ is by no means strongly continuous. Nevertheless, Theorem~\ref
{T1} shows that the limit exists.

A nice by-product of this construction is the following statement
concerning the distribution of $X_t$ for fixed $t$.

\begin{theorem}
\label{T2}
For all $t\geq0$, $x\in\mathbb{L}^{\beta,E}$ and $y\in\mathbb
{L}_\infty^{\beta}$, we have ${\mathbf{E}}
_x[Q_{\langle
X_t,y\rangle}]=Q_{\langle\mathcaligr{S}_tx,y\rangle}$.
In particular, for all $k\in S$, we have
\[
{\mathbf{P}}_x[X_t(k)\in\bolds{\cdot}]=Q_{\mathcaligr{S}_tx(k)}.
\]
\end{theorem}

As an application of Theorem~\ref{T2}, we consider the interface
problem in dimension $d=1$. Assume that $S={\mathbb{Z}}$ and that
${\mathcaligr{A}}
f(k)=\frac12f(k+1)
+\frac12f(k-1)-f(k)$ is the $q$-matrix of symmetric simple random walk
on ${\mathbb{Z}}$. Hence, $a_t$ is the time $t$ transition kernel of
continuous-time rate 1 symmetric simple random walk. Let $u,v>0$ and
assume that $x(k)=(u,0)$ for $k<0$ and $x(k)=(0,v)$ for $k\geq0$. Let
$X$ be the infinite rate mutually catalytic branching process on
${\mathbb{Z}}$
with $X_0=x$. Define
\[
b_{t,1}:=\sup\{k\in{\mathbb{Z}}\dvtx X_{1,t}(k-1)>0\}
\]
and
\[
b_{t,2}:=\inf\{k\in{\mathbb{Z}}\dvtx X_{2,t}(k)>0\}.
\]
We conjecture that $b_{t,1}=b_{t,2}$ almost surely. In this case, the
position $b_t:=b_{t,1}$ could be considered as the interface between
the type 1 population (left) and the type 2 population (right). It is a
challenging task to determine the dynamics of $(b_t)_{t\geq0}$. By work
on the finite branching rate process of \cite{CoxKlenke2000} and \cite
{CoxKlenkePerkins2000}, we should have that $\limsup_{t\to\infty
}b_t=\infty$ and $\liminf_{t\to\infty}b_t=-\infty$. That is,
at any given site, the type changes over and over again at arbitrarily late times.

Theorem~\ref{T2} provides an indication as to what the distribution of
$b_t$ is for fixed $t$.
\begin{corollary}
\label{C1.2}
If $b_{t,1}=b_{t,2}$ almost surely, then
\begin{equation}
\label{E1.14}
{\mathbf{P}}[b_t\leq k]=\frac12+\frac1\pi\arctan\biggl(\frac
{v_t(k)^2-u_t(k)^2}{2u_t(k)v_t(k)}\biggr),
\end{equation}
where
\[
u_t(k):=u\sum_{l=k+1}^\infty a_t(0,l)
\quad\mbox{and}\quad
v_t(k):=v\sum_{l=-k}^\infty a_t(0,l).
\]
In particular, $\operatorname{median}(b_t)\sim\alpha\sqrt{t}$ as $t\to
\infty
$, where $\alpha=\Phi^{-1}(\frac{u}{u+v})$ and $\Phi$ is the
distribution function of the standard normal distribution, and $\lim
_{t\to\infty}{\mathbf{P}}[b_t\leq0]=\frac12+\frac1\pi\arctan
((v^2-u^2)/2uv)$.
\end{corollary}

\begin{pf}
By Theorem~\ref{T2}, we have ${\mathbf{P}}[b_t\leq k]={\mathbf
{P}}[X_{2,t}(k)>0]=Q_{\mathcaligr{S}
_tx(k)}(\{0\}\times(0,\infty))$. By an explicit calculation using the
density of $Q$ (see Lemma~\ref{L2.1.1}), we get \mbox{\rm(\ref
{E1.14})}. The other
two statements follow from the central limit theorem\break for $a_t$.
\end{pf}

\subsection{Outline}
\label{S1.4}
The rest of the paper is organized as follows. In Section~\ref{S2}, we
collect some basic facts about the harmonic measure $Q$ and prove
Proposition~\ref{P1.1}. In Section~\ref{S2b}, we derive a submartingale
related to $X^{\varepsilon }$ and show that the two types of
$X^{\varepsilon }$ are
nonpositively correlated. In Section~\ref{S4}, we show relative
compactness of the family $(X^{\varepsilon }, {\varepsilon }>0)$.
Finally, in Section~\ref
{S5}, we complete the proofs of Theorems~\ref{T1} and \ref{T2}.

\eject
\section{The harmonic measure $Q$}
\label{S2}
\subsection{Harmonic measure and duality}
\label{S2.1}
Recall that $Q_x$ is the harmonic measure for planar Brownian motion in
the upper-right quadrant started at $x\in[0,\infty)^2$ and stopped upon
leaving $(0,\infty)^2$.
If $x=(u,v)\in(0,\infty)^2$, then the harmonic measure $Q_x$ has a
one-dimensional Lebesgue density on $E$
that can be computed explicitly:
\begin{equation}
\label{E2.1.1}
Q_{(u,v)}(d(\bar u,\bar v))=
\cases{\displaystyle \frac4\pi  \frac{\textstyle uv
\bar u}{\textstyle4u^2v^2+(\bar u^2+ v^2-u^2 )^2} \,d\bar u,
&\quad\mbox{if }$\bar v=0$,
\cr
\displaystyle\frac4\pi  \frac
{\textstyle uv \bar v} {\textstyle4u^2v^2+(\bar v^2+u^2-v^2
)^2} \,d\bar v, &\quad\mbox{if }$\bar u=0.$
}
\end{equation}
Furthermore, trivially, we have
$Q_x=\delta_x$ if $x\in E$. Clearly,
\begin{equation}
\label{E2.1.2}
x\mapsto Q_x\qquad \mbox{is continuous. }
\end{equation}

\begin{lemma}
\label{L2.1.1}
For all $u,v>0$ and $c\geq0$, we have
\[
Q_{(u,v)}\bigl(\{0\}\times[c,\infty)\bigr)=\frac12+\frac1\pi\arctan
\biggl(\frac
{v^2-u^2-c^2}{2uv}\biggr).
\]
\end{lemma}

\begin{pf}
This follows from explicitly computing the integral $\int_c^\infty
Q_{(u,v)}(d(0,\break \bar v))$ in \mbox{\rm(\ref{E2.1.1})}.
\end{pf}

Recall $F$ from \mbox{\rm(\ref{E1.9})}. Explicitly computing the
Laplacian with
respect to the first coordinate gives
\[
\biggl(\frac{\partial^2}{(\partial x_1)^2}+\frac{\partial
^2}{(\partial
x_2)^2}\biggr)F(x,y) = 8y_1y_2 F(x,y).
\]
Hence, for $y\in E$, the function $F(\bolds{\cdot},y)$ is harmonic
for planar
Brownian motion $B$ and hence $(F(B_t,y))_{t\geq0}$ is a bounded
martingale. If $\tau$ denotes the first exit time of $B$ from
$(0,\infty
)^2$, then we infer for $x\in[0,\infty)^2$ and $y\in E$ that
\begin{equation}
\label{E2.1.3}
\int F(z,y) Q_x(dz)={\mathbf{E}}_x[F(B_\tau,y)]={\mathbf
{E}}_x[F(B_0,y)]=F(x,y)
\end{equation}
and, similarly (see \cite{KM1}, Corollary 2.3),
\begin{equation}
\label{E2.1.4}
\int F(z,y) Q_x(dz)=\int F(x,z) Q_y(dz)\qquad\mbox{for }x,y\in
[0,\infty)^2.
\end{equation}

Similarly, since linear functions are harmonic for Brownian motion and
using the fact that $p$th moments of $(B_t)_{t\leq\tau}$ are bounded
for $p<2$ (see Lemma~\ref{L2.2.2}), we can derive
\begin{equation}
\label{E2.1.5}
\int z_i  Q_x(dz)=x_i\qquad\mbox{for all }x\in[0,\infty)^2, i=1,2.
\end{equation}
Note that \mbox{\rm(\ref{E2.1.5})} could also be computed explicitly using
Lemma~\ref{L2.2.1}.

Recall that ${\mathbf{Q}}(x,\bolds{\cdot})=\bigotimes_{k\in
S}Q_{x(k)}$ for $x\in
([0,\infty)^2)^S$. From \mbox{\rm(\ref{E2.1.3})}, we
immediately get the
following lemma.
\begin{lemma}
\label{L2.1.2}
For all $x\in\mathbb{L}^{\beta,E}$ and $y\in\mathbb{L}^{f,E}$, we have
\begin{equation}
\label{E2.1.6}
\int H(z,y) {\mathbf{Q}}(x,dz) = H(x,y).
\end{equation}
\end{lemma}

\begin{pf}{Proof of Proposition~\ref{P1.1}}
Note that, due to the definition of $X^{\varepsilon }$ and the chain
rule of
calculus, we have
\[
M^{{\varepsilon },x,y}_t-M^{{\varepsilon },x,y}_s=0\qquad\mbox{for
}s,t\in\bigl[n{\varepsilon },(n+1){\varepsilon }\bigr), n\in{\mathbb{N}}_0.
\]
Hence, the statement of Proposition~\ref{P1.1} is an immediate
consequence of\break Lemma~\ref{L2.1.2}.
\end{pf}

\subsection{Moments of the harmonic measure}
\label{SN2.2}
Since the harmonic measure $Q$ does not possess a second moment, our
proofs will rely on $p$th moment estimates for $p\in(1,2)$. Here, we
collect some of these estimates.
Define $\arctan^\dagger$ as the inverse of the tangent function $\tan
\dvtx [0,\pi]\to\bar{\mathbb{R}}$. That is,
\[
\arctan^\dagger(x)=\arctan(x)+\pi{\mathbh{1}}_{\{x<0\}}.
\]
Note that ${\mathbb{R}}\setminus\{0\}\to[0,\pi]$, $x\mapsto\arctan
^\dagger(1/x)$
can be extended continuously to $x=0$ with the convention that $\arctan
^\dagger(1/0)=\arctan^\dagger(-1/0)=\pi/2$.

\begin{lemma}
\label{L2.2.1}
For all $u,v>0$, we have $\int_Ex_1^2 Q_{(u,v)}(dx)=\infty$ and for
$p\in(0,2)$,
\[
\int_Ex_1^p Q_{(u,v)}(dx)=\frac{(u^2+v^2)^{p/2}
\sin({(p/2)}\arctan^\dagger(
{(2uv/v^2-u^2)}
))}{\sin((\pi/2)p)}.
\]
\end{lemma}

\begin{pf}
This follows from explicitly computing the integral using \mbox{\rm
(\ref{E2.1.1})}.
\end{pf}

\begin{lemma}
\label{L2.2.2}
Let $B=(B_1,B_2)$ be a planar Brownian motion started in $B_0=(u,v)\in
[0,\infty)^2$ and let
\[
\tau=\inf\{t>0\dvtx B_t\notin(0,\infty)^2\}.
\]
Then, for any $p\in(0,2)$, we have
\begin{equation}
\label{E2.2.1}
{\mathbf{E}}[\tau^{p/2}]\leq
\frac{8}{(2-p)(2\pi)^{p/2}}(uv)^{p/2} < \infty.
\end{equation}
Furthermore, for any $p\in(1,2)$, we have
\begin{equation}
\label{E2.2.2}
{\mathbf{E}}[\tau^{p/2}]\leq
\frac{4}{(2\pi)^{p/2}}\frac{p}{(p-1)(2-p)}\min(u^{p-1}v,
uv^{p-1}).
\end{equation}
\end{lemma}

\begin{pf}
By the reflection principle and independence of $B_1$ and $B_2$, we get
\[
{\mathbf{P}}[\tau>t]=4 \mathcaligr{N}_{0,t}(0,u) \mathcaligr{N}_{0,t}(0,v),
\]
where $\mathcaligr{N}_{0,t}(a,b)=(2\pi t)^{-1/2}\int_a^be^{-r^2/2t}\,dr$
is the
centred normal distribution with variance $t$. Hence, for any $p\in(0,2)$,
\begin{eqnarray}
\label{E2.2.3}
{\mathbf{E}}[\tau^{p/2}] &=& \int_0^\infty{\mathbf
{P}}[\tau
>t^{2/p}]\, dt\nonumber
\\[-8pt]\\[-8pt]
&\leq&  4\int_0^\infty
\bigl(1\wedge u(2\pi)^{-1/2} t^{-1/p}\bigr)\bigl(1\wedge v(2\pi
)^{-1/2}t^{-1/p}\bigr)\, dt.\nonumber
\end{eqnarray}
We can continue this inequality as
\[
\leq\frac{4}{(2\pi)^{p/2}}(uv)^{p/2}+\frac{2uv}{\pi}\int
_{(uv/2\pi
)^{p/2}}^\infty t^{-2/p}\, dt = \frac{8}{(2-p)(2\pi)^{p/2}}(uv)^{p/2}.
\]
This gives \mbox{\rm(\ref{E2.2.1})}. For $p\in(1,2)$, we can
continue \mbox{\rm(\ref{E2.2.3})} as
\begin{eqnarray}
\label{E2.2.4}
&\leq&4u(2\pi)^{-1/2}\int_0^{v^p/(2\pi)^{p/2}}t^{-1/p}\, dt + \frac
{2uv}{\pi}\int_{v^p/(2\pi)^{p/2}}^\infty t^{-2/p}\, dt\nonumber
\\[-8pt]\\[-8pt]
& =& \frac{4}{(2\pi)^{p/2}} \frac{p}{(p-1)(2-p)} uv^{p-1}.\nonumber
\end{eqnarray}
Interchanging the roles of $u$ and $v$ in \mbox{\rm(\ref{E2.2.4})}
gives \mbox{\rm(\ref{E2.2.2})}.
\end{pf}

\begin{lemma}
\label{L2.2.3}
For $p\in(1,2)$, there exists a constant $C_p<\infty$ such that for
every $x\in[0,\infty)^2$ and $i=1,2$, we have
\begin{equation}
\label{E2.2.5}
\int_E y_i^p Q_x(dy)\geq x_i^p
\end{equation}
and
\begin{equation}
\label{E2.2.6}
\int_E|y_i-x_i|^p Q_x(dy) \leq  C_p   \min(x_1^{p-1}x_2,
x_1x_2^{p-1}) \leq  C_p  (x_1x_2)^{p/2}.
\end{equation}
\end{lemma}

\begin{pf}
By the Burkholder--Davis--Gundy inequality (see, e.g.,
\cite{DellacherieMeyer1983},
Theorem~VII.92) and Lemma~\ref{L2.2.2},
$(B_{i,t})_{t\leq\tau}$ is a uniformly integrable martingale. Hence,
by Jensen's inequality,
\[
x_i^p = {\mathbf{E}}_x[B_{i,t}]^p \leq {\mathbf{E}}_x
[B_{i,t}^p] = \int
_Ey_i^p Q_x(dy).
\]

The claim \mbox{\rm(\ref{E2.2.6})} could be checked either by a
direct computation
using Lemma~\ref{L2.2.1} or by proceeding as follows. Let $B$ and
$\tau
$ be as in Lemma~\ref{L2.2.2}. Using the Burkholder--Davis--Gundy
inequality and then Lemma~\ref{L2.2.2}, we get
\begin{eqnarray*}
\int_E|y_i-x_i|^p Q_x(dy) &=& {\mathbf{E}}_x[|B_{i,\tau
}-x_i|^p]\leq
(4p)^p {\mathbf{E}}_x[\tau^{p/2}]
\\
&\leq&\frac
{(4p)^{p+1}}{(p-1)(2-p)(2\pi)^{p/2}} \min
(x_1^{p-1}x_2,x_1x_2^{p-1}).
\end{eqnarray*}
\upqed\end{pf}

\section{The approximating process $X^{\varepsilon }$}
\label{S2b}
\subsection{Martingale property of $X^{\varepsilon }$}
\begin{proposition}
\label{P3.1.1}
Let $x\in\mathbb{L}^{\beta,E}$ and $k\in S$.
Define the process $N^{{\varepsilon },x}$ for $i=1,2$, $k\in S$ and
$t\geq0$ by
\[
N^{{\varepsilon },x}_{i,t}(k):=X^{\varepsilon
}_{i,t}(k)-X^{\varepsilon }_{i,0}(k)-\int_0^t({\mathcaligr{A}}
X^{\varepsilon }
_{i,s})(k) \,ds.
\]
\begin{longlist}[(ii)]
\item[(i)]
For each $i=1,2$ and $k\in S$, the process $(N^{{\varepsilon }
,x}_{i,t}(k))_{t\geq
0}$ is a martingale with respect to the natural filtration. In particular,
\begin{equation}
\label{E3.1.1}
{\mathbf{E}}_x[X^{\varepsilon }_{i,t}(k)]=(\mathcaligr{S}_tx_i)(k)\qquad
\mbox{for all }t\geq0, k\in S, i=1,2.
\end{equation}
\item[(ii)]
Define $\lambda:=\sup_{k\in S}(-{\mathcaligr{A}}(k,k))$ and note that
$|\lambda
|<\infty$ by assumption \mbox{\rm(\ref{E1.3})}. Define
\[
Z^{\varepsilon }_{i,t}(k):=e^{-{\mathcaligr{A}}(k,k)t}X^{\varepsilon }_{i,t}(k)
\]
and
\[
\bar Z^{\varepsilon }_{i,t}:=e^{\lambda t}\| X^{\varepsilon
}_{i,t}\|_\beta.
\]
$Z^{\varepsilon }_{i}(k)$ and $\bar Z^{\varepsilon }_{i}$ are then
nonnegative submartingales.
\end{longlist}
\end{proposition}

\begin{pf}
{(i)} This is an immediate consequence of the definition
of $X^{\varepsilon }$ and \mbox{\rm(\ref{E2.1.5})}.

{(ii)} Since ${\mathcaligr{A}}(k,l)\geq0$ for all $k\neq
l$, we have
\[
\frac{d}{dt} Z^{\varepsilon }_{i,t}(k)=\sum_{l\neq k}{\mathcaligr
{A}}(k,l)Z^{\varepsilon }
_{i,t}(l)\geq
0\qquad\mbox{for }t\in\bigl(n{\varepsilon },(n+1){\varepsilon }\bigr).
\]
Together with \mbox{\rm(\ref{E2.1.5})}, this shows that
$Z^{\varepsilon }_{i}$ is a submartingale.
As a sum of submartingales, $\bar Z^{\varepsilon }_{i}$ is also a
submartingale.
\end{pf}

\begin{corollary}
\label{C3.1.2}
For every $K,T>0$ and any set $G\subset S$, we have
\begin{equation}
\label{E3.1.2}
\quad{\mathbf{P}}_x\Bigl[\sup_{t\in[0,T]}
\| (X^{\varepsilon }_{1,t}+X^{\varepsilon }_{2,t}) {\mathbh
{1}}_G\|_\beta\geq K\Bigr]\leq
K^{-1}e^{\lambda T} \bigl\|\bigl(\mathcaligr{S}_T(x_1+x_2)\bigr)
{\mathbh{1}}_G\bigr\|_\beta.
\end{equation}
In particular,
\begin{equation}
\label{E3.1.3}
{\mathbf{P}}_x\Bigl[\sup_{t\in[0,T]}
\| (X^{\varepsilon }_{1,t}+X^{\varepsilon }_{2,t})\|_\beta
\geq K\Bigr]\leq
K^{-1}e^{(\lambda+M)T}\|x_1+x_2\|_\beta.
\end{equation}
\end{corollary}

\begin{pf}
This is an immediate consequence of Proposition~\ref{P3.1.1} and Doob's
inequality.
\end{pf}

\subsection{One-dimensional distributions}
\label{S3.2}
\begin{lemma}
\label{L3.2.1}
Let $a(1),a(2),\ldots$ be nonnegative numbers and let
$x(1),x(2),\break\ldots
\in[0,\infty)^2$ be such that
\[
\bar x:=\langle a,x\rangle=\sum_{k=1}^\infty a(k) x(k)\in[0,\infty)^2.
\]
Let $\xi(1),\xi(2),\ldots$ be independent random variables with
${\mathbf{P}}
[\xi
(k)\in\bolds{\cdot}]=Q_{x(k)}$. Define $\bar\xi:=\langle a,\xi
\rangle=\sum
_{k=1}^\infty a(k)\xi(k)$ and assume that $X$ is an $E$-valued random
variable such that ${\mathbf{P}}[X\in\bolds{\cdot}\hspace
{0.8pt}|\hspace{0.8pt}
\bar\xi ]=Q_{\bar\xi
}\raisebox{2pt}{\rule{0pt}{8pt}}$. Then ${\mathbf{P}}[X\in\bolds
{\cdot}]=Q_{\bar x}$.
In other words, ${\mathbf{E}}[Q_{\bar\xi}]=Q_{\bar x}$.
\end{lemma}

\begin{pf}
First, note that ${\mathbf{E}}[\xi_i(k)]=x_i(k)$ and, hence, $\bar
\xi\in
[0,\infty
)^2$ almost surely. Recall $F$ from \mbox{\rm(\ref{E1.9})}. By \mbox
{\rm(\ref{E2.1.3})}, for
all $y\in E$, we have
\begin{eqnarray*}
{\mathbf{E}}[F(X,y)]&=&{\mathbf{E}}[F(\bar\xi,y)]=\prod
_{k=1}^\infty{\mathbf{E}}[F(\xi(k),a(k)y)]
=\prod_{k=1}^\infty F(x(k),a(k)y)
\\
&=&F(\bar x,y)=\int_E F(z,y)
Q_{\bar x}(dz).
\end{eqnarray*}
Since $F(\bolds{\cdot},y)$, $y\in E$, is measure-determining (see
\cite{KM1},
Corollary 2.4), this yields the claim.
\end{pf}

\begin{corollary}
\label{C3.2.2}
For any ${\varepsilon }>0$, $n\in{\mathbb{N}}_0$ and $k\in S$, we have
\[
{\mathbf{P}}[X^{\varepsilon }_{n{\varepsilon }}(k)\in\bolds
{\cdot}]=Q_{\mathcaligr{S}_{n{\varepsilon }}x(k)}.
\]
\end{corollary}

\begin{pf}
Fix $n\in{\mathbb{N}}$.
We show by induction on $m$ that
\[
{\mathbf{P}}[X^{\varepsilon }_{n{\varepsilon }}(k)\in\bolds
{\cdot}\hspace{0.8pt}|\hspace{0.8pt}X^{\varepsilon }
_{m{\varepsilon }}]
=Q_{(\mathcaligr{S}_{(n-m){\varepsilon }} X^{\varepsilon
}_{m{\varepsilon }})(k)}\qquad\mbox{for all }m=0,\ldots,n.
\]
For the induction base $m=n$, this is true by the definition of
$X^{\varepsilon }
$. Now, assume that we have shown the statement for some $m\geq1$.
Using the induction hypothesis in the first line and Lemma~\ref{L3.2.1}
in the second line, we get
\begin{eqnarray*}
{\mathbf{P}}\bigl[X^{\varepsilon }_{n{\varepsilon }}(k)\in\bolds
{\cdot}|X^{\varepsilon }
_{(m-1){\varepsilon }}\bigr]
&=&{\mathbf{E}}\bigl[Q_{(\mathcaligr{S}_{(n-m){\varepsilon }}
X^{\varepsilon }_{m{\varepsilon }})(k)}|X^{\varepsilon }
_{(m-1){\varepsilon }}
\bigr]
\\
&=&{\mathbf{E}}\bigl[Q_{(\mathcaligr{S}_{(n-m){\varepsilon }}
X^{\varepsilon }_{m{\varepsilon }-})(k)}|X^{\varepsilon }
_{(m-1){\varepsilon }
}\bigr]
\\
&=&Q_{(\mathcaligr{S}_{(n-(m-1)){\varepsilon }} X^{\varepsilon
}_{(m-1){\varepsilon }})(k)}.
\end{eqnarray*}
Note that we have used the fact that $X^{\varepsilon }_{m{\varepsilon
}-}=\mathcaligr{S}_{\varepsilon }X^{\varepsilon }
_{(m-1){\varepsilon }}$ in the last line.
\end{pf}

\begin{corollary}
\label{C3.2.3}
Let $y\in\mathbb{L}_\infty^{\beta}$.
Then ${\mathbf{E}}_x[Q_{\langle X^{\varepsilon }_t,y\rangle
}]=Q_{\langle\mathcaligr{S}
_tx(k),y\rangle}$.
\end{corollary}

\begin{pf}
The proof is similar to the proof of Corollary~\ref{C3.2.2}. (Note that
$\langle X^{\varepsilon }_t,y\rangle\in[0,\infty)^2$ almost surely since
$X^{\varepsilon }
_t\in\mathbb{L}^{\beta,2}$ almost surely.)
\end{pf}

\subsection{Correlations}
\begin{lemma}
\label{L3.3.1}
Let $Y$ and $Z$ be nonpositively correlated nonnegative random
variables and assume that $h\dvtx[0,\infty)\to[0,\infty)$ is concave
and monotone increasing. Then
$
{\mathbf{E}}[Yh(Z)]\leq{\mathbf{E}}[Y] h({\mathbf{E}}[Z])
$.
\end{lemma}

\begin{pf}
If ${\mathbf{E}}[Z]=0$, then we even have equality. Now, assume that
${\mathbf{E}}[Z]>0$.
By concavity of $h$, there exists a $b\in{\mathbb{R}}$ such that for
all $z\geq0$,
\[
h(z)\leq h({\mathbf{E}}[Z])+(z-{\mathbf{E}}[Z]) b.
\]
Since $h$ is nondecreasing, we have $b\geq0$ and thus
\[
{\mathbf{E}}[Y h(Z)]
 \leq  {\mathbf{E}}\bigl[Y \bigl(h({\mathbf{E}}[Z])+(Z-{\mathbf
{E}}[Z]) b\bigr)\bigr]
 \leq  {\mathbf{E}}[Y] h({\mathbf{E}}[Z]).
\]
\upqed\end{pf}

\begin{lemma}
\label{L3.3.2}
For any ${\varepsilon }>0$, $n\in{\mathbb{N}}_0$ and $k\in S$, the
random variables
$X^{\varepsilon }
_{1,n{\varepsilon }}(k)$ and $X^{\varepsilon }_{2,n{\varepsilon
}}(k)$ are nonpositively correlated, in
the sense that
\begin{eqnarray}
\label{E3.3.1}
{\mathbf{E}}_x[X^{\varepsilon }_{1,n{\varepsilon
}}(k)X^{\varepsilon }_{2,n{\varepsilon }}(k)] &\leq&
{\mathbf{E}}_x[X^{\varepsilon }_{1,n{\varepsilon }}(k)]
{\mathbf{E}}_x[X^{\varepsilon }_{2,n{\varepsilon }
}(k)]\nonumber
\\[-8pt]\\[-8pt]
&=&
(\mathcaligr{S}_{n{\varepsilon }}x_1(k))(\mathcaligr{S}_{n{\varepsilon }}x_2(k)).\nonumber
\end{eqnarray}
\end{lemma}

\begin{pf}
Let $t\geq0$. Recall that ${\mathcaligr{F}}$ is the natural filtration
of $X^{\varepsilon }
$. Then
\begin{eqnarray*}
&&{\mathbf{E}}_x\bigl[\mathcaligr{S}_tX^{\varepsilon }_{1,n{\varepsilon
}}(k) \mathcaligr{S}_tX^{\varepsilon }_{2,n{\varepsilon }}(k)
|
{\mathcaligr{F}}
_{(n-1){\varepsilon }}\bigr]
\\
&&\qquad
=\sum_{l_1\neq l_2}a_t(k,l_1) a_t(k,l_2) {\mathbf{E}}_x
\bigl[X^{\varepsilon }_{1,n{\varepsilon }
}(l_1) X^{\varepsilon }_{2,n{\varepsilon }}(l_2)
|{\mathcaligr{F}}_{(n-1){\varepsilon }}\bigr]
\\
&&\qquad =\sum_{l_1\neq l_2}a_t(k,l_1) a_t(k,l_2) \bigl(\mathcaligr
{S}_{\varepsilon }X^{\varepsilon }
_{1,(n-1){\varepsilon }
}\bigr)(l_1) \bigl(\mathcaligr{S}_{\varepsilon }X^{\varepsilon
}_{2,(n-1){\varepsilon }}\bigr)(l_2)
\\
&&\qquad \leq\sum_{l_1, l_2}a_t(k,l_1)\bigl(\mathcaligr{S}_{\varepsilon
}X^{\varepsilon }_{1,(n-1){\varepsilon }}\bigr)(l_1)
a_t(k,l_2)\bigl(\mathcaligr{S}_{\varepsilon }X^{\varepsilon
}_{2,(n-1){\varepsilon }}\bigr)(l_2)
\\
&&\qquad =\mathcaligr{S}_{t+{\varepsilon }} X^{\varepsilon
}_{1,(n-1){\varepsilon }}(k) \mathcaligr{S}_{t+{\varepsilon }}
X^{\varepsilon }
_{2,(n-1){\varepsilon }}(k).
\end{eqnarray*}
Inductively, we get
\[
{\mathbf{E}}_x[\mathcaligr{S}_tX^{\varepsilon }_{1,n{\varepsilon
}}(k) \mathcaligr{S}_tX^{\varepsilon }_{2,n{\varepsilon }}(k)]
\leq\mathcaligr{S}_{t+n{\varepsilon }}x_1(k) \mathcaligr
{S}_{t+n{\varepsilon }}x_2(k).
\]
Applying this with $t=0$ yields the claim.
\end{pf}

\section{Tightness}
\label{S4}
The goal of this section is to show the following proposition.
\begin{proposition}
\label{P4.1}
The family of processes $(X^{\varepsilon })_{{\varepsilon }>0}$ is
relatively compact in the Skorohod spaces of c{\`a}dl{\`a}g functions
$D([0,\infty);\mathbb{L}^{\beta,2})$.
\end{proposition}

By Prohorov's theorem, in order to show relative compactness of
$(X^{\varepsilon }
)$, it is enough to show tightness of $(X^{\varepsilon })$.

The strategy of proof is to check the compact containment condition for
$X^{\varepsilon }$ (Lemma~\ref{L4.4}) and then use Aldous's tightness
criterion
for functions $h(X^{\varepsilon }_t)$, where $h\dvtx \mathbb{L}^{\beta
,2}\to{\mathbb{R}}$ is Lipschitz continuous
and depends on only finitely many coordinates.

We start by collecting some basic facts about compact sets and
separating function spaces. The proofs of the following statements are
standard and are therefore omitted here.
\begin{lemma}
\label{L4.2} A set $C\subset\mathbb{L}^{\beta,2}$ is relatively
compact if and only if
the following hold:
\begin{longlist}[(ii)]
\item[(i)]
$B_C:=\sup_{x\in C}\|x_1+x_2\|_\beta<\infty$;
\item[(ii)]
for any $\eta>0$, there exists a finite subset $S_\eta\subset S$ such
that $\sup_{x\in C}\|(x_1+x_2) {\mathbh{1}}_{S\setminus S_\eta
}\|
_\beta
<\eta$.
\end{longlist}
\end{lemma}

\begin{lemma}
\label{L4.3}Let $C_b(\mathbb{L}^{\beta,2};{\mathbb{R}})$ be the
space of real-valued bounded
continuous functions $\mathbb{L}^{\beta,2}\to{\mathbb{R}}$ with the
topology of uniform
convergence on compact sets.
Denote by $\mathrm{Lip}_f(\mathbb{L}^{\beta,2};{\mathbb{R}})$ the
space of Lipschitz continuous
bounded functions $\mathbb{L}^{\beta,2}\to{\mathbb{R}}$ that depend
on only finitely many
coordinates. Then $\mathrm{Lip}_f(\mathbb{L}^{\beta,2};{\mathbb
{R}})\subset C_b(\mathbb{L}^{\beta,2};{\mathbb{R}})$ is dense.
\end{lemma}

\begin{lemma}[(Compact containment condition)]
\label{L4.4}
Fix $x\in\mathbb{L}^{\beta,2}$. For any $\eta>0$ and $T>0$, there
exists a compact set
$\Gamma\subset\mathbb{L}^{\beta,2}$ such that
\begin{equation}
\label{E4.1}
{\mathbf{P}}_x\bigl[X^{\varepsilon }_t\in\Gamma\mbox{ for all
}t\in[0,T]\bigr]\geq1-\eta
\qquad\mbox{for all }{\varepsilon }>0.
\end{equation}
\end{lemma}

\begin{pf}
Let $T>0$ and $\eta>0$. Recall $M$ from \mbox{\rm(\ref{E1.4})} and
$\lambda$ from
Proposition~\ref{P3.1.1}(ii). Choose a $K>\frac{2}{\eta}
e^{(\lambda
+M)T} \|x_1+x_2\|_\beta$ and let $A_K:=\{y\in\mathbb{L}^{\beta
,2}\dvtx \|y_1+y_2\|
_\beta
<K\}$. According to Corollary~\ref{C3.1.2}, we have
\[
{\mathbf{P}}_x\bigl[X^{\varepsilon }_t\in A_K\mbox{ for all }t\in
[0,T]\bigr]\geq1-\frac
{\eta}{2}.
\]
Now, for any $n\in{\mathbb{N}}$, choose a finite $S_n\subset S$ such that
\[
n e^{\lambda T}\|\mathcaligr{S}_T(x_1+x_2) {\mathbh
{1}}_{S\setminus S_n}\|
_\beta
< 2^{-n-1} \eta
\]
and define
\[
B_n:=\{y\in\mathbb{L}^{\beta,2}\dvtx \|(y_1+y_2){\mathbh
{1}}_{S\setminus S_n}\|_\beta
<1/n\}.
\]
According to Corollary~\ref{C3.1.2}, we have
\[
{\mathbf{P}}_x\bigl[X^{\varepsilon }_t\in B_n\mbox{ for all }t\in
[0,T]\bigr] \geq
1-2^{-n-1}
\eta.
\]
Now, let $\Gamma$ by the closure of $A_K\cap\bigcap_{n=1}^\infty B_n$.
Then
\[
{\mathbf{P}}_x\bigl[X^{\varepsilon }_t\in\Gamma\mbox{ for all
}t\in[0,T]\bigr] \geq
1-\eta
\]
and, by Lemma~\ref{L4.2}, $\Gamma$ is compact.
\end{pf}

\begin{lemma}
\label{L4.5}
Fix $h\in\mathrm{Lip}_f(\mathbb{L}^{\beta,2};{\mathbb{R}})$. For
${\varepsilon }>0$, define the process
$Y^{\varepsilon }$ by
\[
Y^{\varepsilon }_t:=h(X^{\varepsilon }_t),\qquad t\geq0.
\]
$(Y^{\varepsilon })_{{\varepsilon }>0}$ is then tight in the
Skorohod space $D([0,\infty
);{\mathbb{R}}
)$ of c{\`a}dl{\`a}g functions $[0,\infty)\to{\mathbb{R}}$.
\end{lemma}

\begin{pf}
The idea is to use Aldous's criterion for tightness in $D([0,\infty
);{\mathbb{R}})$.
As~$h$ is bounded, $(Y^{\varepsilon }_t)_{{\varepsilon }>0}$ is tight
for each $t\geq0$. Hence,
by Aldous's criterion (see, e.g., \cite{Aldous1978}, equation (13), or
\cite{JacodShiryaev1987}, Section VI.4a), we need to show the
following: for any $\eta>0$ and $T>0$, there exist $\delta>0$ and
${\varepsilon }
_0>0$ such that, for any stopping time $\tau\leq T$, we have
\begin{equation}
\label{E4.2}
\sup_{\delta'\in[0,\delta]}
 \sup_{{\varepsilon }\in(0,{\varepsilon }_0]}
{\mathbf{P}}_x[|Y^{\varepsilon }_{\tau+\delta
'}-Y^{\varepsilon }_\tau|>\eta]\leq\eta.
\end{equation}
Since $h$ is Lipschitz continuous and depends on only finitely many
coordinates, it is enough to consider the case where $h(x)=x_i(k)$ for
some $k\in S$ and $i=1,2$. Using Markov's inequality, it is enough to
show that for any $\eta>0$ and $T>0$, there exist $\delta>0$ and
${\varepsilon }
_0>0$ such that for any stopping time $\tau\leq T$, we have
\begin{equation}
\label{E4.3}
\sup_{\delta'\in[0,\delta]}
 \sup_{{\varepsilon }\in(0,{\varepsilon }_0]}
{\mathbf{E}}_x[|X^{\varepsilon }_{i,\tau+\delta
'}(k)-X^{\varepsilon }_{i,\tau}(k)|]\leq
\eta.
\end{equation}
Define
\[
N:=\lfloor\tau/{\varepsilon }\rfloor\quad\mbox{and}\quad
N':=\lfloor(\tau+\delta')/{\varepsilon }
\rfloor.
\]
Then
\[
{\mathbf{E}}_x[|X^{\varepsilon }_{i,\tau+\delta
'}(k)-X^{\varepsilon }_{i,\tau}(k)|]
\leq E_1+E_2+E_3+E_4,
\]
where
\begin{eqnarray*}
E_1&:=&{\mathbf{E}}_x[|X^{\varepsilon }_{i,\tau
}(k)-X^{\varepsilon }_{i,N{\varepsilon }}(k)|],\\
E_2&:=&{\mathbf{E}}_x[|X^{\varepsilon }_{i,\tau+\delta
'}(k)-X^{\varepsilon }_{i,N'{\varepsilon }
}(k)|],\\
E_3&:=&{\mathbf{E}}_x[|X^{\varepsilon }_{i,N'{\varepsilon
}}(k)-X^{\varepsilon }_{i,N{\varepsilon }}(k)|].
\end{eqnarray*}
Now, by \mbox{\rm(\ref{E1.7})}, we get
\begin{eqnarray*}
E_1&=&{\mathbf{E}}_x[|\mathcaligr{S}_{\tau-N{\varepsilon
}}X^{\varepsilon }_{i,N{\varepsilon }}(k)-X^{\varepsilon
}_{i,N{\varepsilon }
}(k)|]
\\
&\leq&{\mathbf{E}}_x\biggl[\int_0^{\tau-N{\varepsilon }}|{\mathcaligr
{A}}\mathcaligr{S}_{s}X_{i,N{\varepsilon }}(k)|
\,ds\biggr]
\\
&\leq&\frac{Me^{\delta M}\delta}{\beta(k)} {\mathbf{E}}_x[\|
X^{\varepsilon }
_{i,N{\varepsilon }}\|
_\beta]
\leq\frac{Me^{(T+2\delta)M}}{\beta(k)} \|x_i\|_\beta  \delta.
\end{eqnarray*}
Similarly, we get
\[
E_2
\leq\frac{Me^{(T+2\delta)M}}{\beta(k)} \|x_i\|_\beta \delta.
\]
Note that $N'-N$ takes only the values $\lfloor\delta'/{\varepsilon
}\rfloor$ and
$\lceil\delta'/{\varepsilon }\rceil$. Hence, $E_3\leq E_3'+E_3''$, where
\[
E_3':={\mathbf{E}}_x\bigl[\bigl|X^{\varepsilon }_{i,(N+\lfloor\delta
'/{\varepsilon }\rfloor){\varepsilon }
}(k)-X^{\varepsilon }
_{i,N{\varepsilon }}(k)\bigr|\bigr]
\]
and
\[
E_3'':={\mathbf{E}}_x\bigl[\bigl|X^{\varepsilon }_{i,(N+\lceil\delta
'/{\varepsilon }\rceil){\varepsilon }
}(k)-X^{\varepsilon }
_{i,N{\varepsilon }}(k)\bigr|\bigr].
\]
Define
\[
\bar E_3'':={\mathbf{E}}_x\bigl[\bigl|X^{\varepsilon }_{i,(N+\lceil
\delta'/{\varepsilon }\rceil
){\varepsilon }
}(k)-\mathcaligr{S}_{\lceil\delta'/{\varepsilon }\rceil{\varepsilon
}}X^{\varepsilon }_{i,N{\varepsilon }}(k)
\bigr|\bigr].
\]
Using the triangle inequality and proceeding as for $E_1$, we get
\[
E_3''\leq\bar E_3''+{\mathbf{E}}_x\bigl[\bigl|X^{\varepsilon
}_{i,N{\varepsilon }}(k)-\mathcaligr{S}_{\lceil
\delta
'/{\varepsilon }\rceil{\varepsilon }}X^{\varepsilon
}_{i,N{\varepsilon }}(k)\bigr|\bigr]\leq\bar E_3+\frac
{Me^{(T+2\delta)M}}{\beta(k)} \|x_i\|_\beta  \delta.
\]
Fix a $p\in(1,2)$. Using the Markov property of $X^{\varepsilon }$ and
conditioning on $X_{N{\varepsilon }}^{\varepsilon }$, by
Corollary~\ref{C3.2.2}
and Jensen's inequality, we get
\begin{eqnarray*}
\bar E_3''
&=&{\mathbf{E}}_x\biggl[\int_E |y_i-X^{\varepsilon }_{i,N{\varepsilon
}}(k)| Q_{\mathcaligr{S}_{\lceil\delta
'/{\varepsilon }
\rceil{\varepsilon }}X^{\varepsilon }_{N{\varepsilon }}(k)}(dy)
\biggr]
\\
&\leq&\biggl({\mathbf{E}}_x\biggl[\int_E |y_i-X^{\varepsilon
}_{i,N{\varepsilon }}(k)|^p Q_{\mathcaligr{S}
_{\lceil
\delta'/{\varepsilon }\rceil{\varepsilon }}X^{\varepsilon
}_{N{\varepsilon }}(k)}(dy)\biggr]\biggr)^{1/p}.
\end{eqnarray*}
Applying Lemma~\ref{L2.2.3}, there exists a constant $C=C_p<\infty$
such that
\begin{eqnarray*}
(\bar E_3'')^p &\leq&  C
{\mathbf{E}}_x\bigl[
\bigl(\mathcaligr{S}_{\lceil\delta'/{\varepsilon }\rceil{\varepsilon
}}X^{\varepsilon }_{1,N{\varepsilon }}(k)
\bigr)^{p-1}\mathcaligr{S}
_{\lceil\delta'/{\varepsilon }\rceil{\varepsilon }}X^{\varepsilon
}_{2,N{\varepsilon }}(k) {\mathbh{1}}_{\{X^{\varepsilon }
_{2,N{\varepsilon }
}(k)=0\}}\bigr]
\\
&&{} +C {\mathbf{E}}_x\bigl[
\bigl(\mathcaligr{S}_{\lceil\delta'/{\varepsilon }\rceil{\varepsilon
}}X^{\varepsilon }_{2,N{\varepsilon }}(k)
\bigr)^{p-1}\mathcaligr{S}
_{\lceil\delta'/{\varepsilon }\rceil{\varepsilon }}X^{\varepsilon
}_{1,N{\varepsilon }}(k) {\mathbh{1}}_{\{X^{\varepsilon }
_{1,N{\varepsilon }
}(k)=0\}}\bigr].
\end{eqnarray*}
By symmetry, it is enough to consider the first summand. Since the
first and the second type are nonpositively correlated (Lemma~\ref
{L3.3.2}), by Lemma~\ref{L3.3.1} [with $h(z)=z^{p-1}$], the first
summand can be estimated by
\begin{eqnarray*}
&&{\mathbf{E}}_x\bigl[
\bigl(\mathcaligr{S}_{\lceil\delta'/{\varepsilon }\rceil{\varepsilon
}}X^{\varepsilon }_{1,N{\varepsilon }}(k)
\bigr)^{p-1}M\delta e^{M\delta}\beta(k)^{-1}\|X^{\varepsilon
}_{N,{\varepsilon },2}\|_\beta
\bigr]
\\
&&\qquad\leq{\mathbf{E}}_x\bigl[
\bigl(\mathcaligr{S}_{\lceil\delta'/{\varepsilon }\rceil{\varepsilon
}}X^{\varepsilon }_{1,N{\varepsilon }}(k)
\bigr)^{p-1}\bigr]
 M e^{M\delta}\delta\beta(k)^{-1}{\mathbf{E}}_x[\|
X^{\varepsilon }_{N,{\varepsilon },2}\|
_\beta
]
\\
&&\qquad \leq\bigl(e^{M(T+\delta+{\varepsilon }_0)}\|x_1\|_\beta\bigr)^{p-1}
 M\delta\beta(k)^{-1}e^{M(T+2\delta)}\|x_2\|_\beta.
\end{eqnarray*}
The estimate for $E_3'$ is analogous.
Summing up, by choosing $\delta$ sufficiently small (independently of
${\varepsilon }\leq{\varepsilon }_0$), we can get $E_j<\eta/3$,
$j=1,2,3$ and hence \mbox{\rm(\ref{E4.3})}.
\end{pf}

\begin{pf*}{Proof of Proposition~\ref{P4.1}}
The space $\mathbb{L}^{\beta,2}$ is Polish and hence so is the
Skorohod space
$D([0,\infty);\mathbb{L}^{\beta,2})$ of c{\`a}dl{\`a}g paths
$[0,\infty)\to\mathbb{L}^{\beta,2}$ (see
\cite{EthierKurtz1986}, Chapter~III.5). Hence, by Prohorov's theorem,
it is enough to show tightness of $(X^{\varepsilon })_{{\varepsilon
}>0}$ in $D([0,\infty
);\mathbb{L}^{\beta,2})$. By \cite{EthierKurtz1986}, Theorem
III.9.1, it is enough to
check two conditions:
\begin{longlist}[(ii)]
\item[(i)] the compact containment condition---this is done in
Lemma~\ref{L4.4};
\item[(ii)] there is a dense (in the topology of uniform convergence on
compacts) space $H\subset C_b(\mathbb{L}^{\beta,2};{\mathbb{R}})$
such that for every $h\in H$,
the family $h(X^{\varepsilon })$, ${\varepsilon }>0$, is tight in
$D([0,\infty);{\mathbb{R}})$---we
have checked this for $H=\operatorname{Lip}_f(\mathbb{L}^{\beta
,2};{\mathbb{R}})$ in Lemmas~\ref{L4.3}
and~\ref{L4.5}.
\end{longlist}
\upqed\end{pf*}

\section{The martingale problem}
\label{S5}
In this section, we complete the proofs of Theorems~\ref{T1} and
\ref{T2}.

\subsection[Proof of Theorem 1]{Proof of Theorem~\protect\ref{T1}}
\label{S5.1}
From Proposition~\ref{P4.1} we know that $X^{\varepsilon }$,
${\varepsilon }>0$, is weakly
relatively compact. From Theorem~\ref{T0}, we know that the martingale
problem \mbox{\rm(\ref{MP1})} has a unique solution. Hence, it
remains to show that
any weak limit point of~$X^{\varepsilon }$, ${\varepsilon }>0$, is a
solution of \mbox{\rm(\ref{MP1})}.

Let $x\in\mathbb{L}^{\beta,E}$. Fix a sequence ${\varepsilon
}_n\downarrow0$ such that $X^{{\varepsilon }_n}$
converges and denote the limit by $X$. Without loss of generality, we
may assume that the processes are defined on one probability space such
that $X^{{\varepsilon }_n}\stackrel{n \rightarrow\infty
}{\longrightarrow}X$ almost surely. Let $y\in\mathbb{L}^{f,E}$ and define
$M^{x,y}$ as in \mbox{\rm(\ref{MP1})} and $M^{{\varepsilon },x,y}$
as in \mbox{\rm(\ref{E1.13})}. We know
from Proposition~\ref{P1.1} that $M^{{\varepsilon }_n,x,y}$ is a martingale.
Hence, it is enough to show that
\begin{equation}
\label{E5.1} M^{{\varepsilon }_n,x,y}_t\stackrel{n \rightarrow
\infty}{\longrightarrow}M^{x,y}_t \qquad\mbox{in }L^1\mbox{ for all
}t\geq0.
\end{equation}
Note that the integrand in \mbox{\rm(\ref{E1.13})} converges
pointwise to the
integrand in \mbox{\rm(\ref{MP1})}. Since $H$ is bounded, in order
to show \mbox{\rm(\ref{E5.1})}, it is enough to show that $\langle
\hspace*{-0.13em}\langle
{\mathcaligr{A}}X^{{\varepsilon }_n}_s,y\rangle\hspace
*{-0.13em}\rangle $ is
uniformly integrable (with respect to Lebesgue measure on $[0,t]$ and
${\mathbf{P}}_x$). Let $p\in(1,2)$. Since $y(k)\neq0$ for only
finitely many
$k\in S$, it is enough to show that for $i=1,2$ and $t>0$, we have
\begin{equation}
\label{E5.2}
\sup_{{\varepsilon }>0}\sup_{s\in[0,t]}{\mathbf{E}}
[|{\mathcaligr{A}}X^{\varepsilon }_{i,s}(k)|^p
]<\infty.
\end{equation}
Recall that $|{\mathcaligr{A}}X^{\varepsilon }_{i,s}(k)|\leq M\|
X^{\varepsilon }_{i,s}\|_\beta/\beta(k)$.
Let $Z$ be an $E$-valued random variable such that ${\mathbf{P}}[Z\in
\bolds{\cdot}
\hspace{0.8pt}|\hspace{0.8pt}
X^{\varepsilon }]=Q_{\|X^{\varepsilon }_s\|_\beta}$. Then
${\mathbf{E}}[Z_i^p]\geq{\mathbf{E}}_x[\|X^{\varepsilon }_{i,s}\|
_\beta^p]$, by Lemma~\ref{L2.2.3}.
However, by Corollary~\ref{C3.2.3}, we have ${\mathbf{P}}_x[Z\in
\bolds{\cdot}]=Q_{\|
\mathcaligr{S}
_sx\|_\beta}$. Hence, again by Lemma~\ref{L2.2.3} and using \mbox
{\rm(\ref{E1.7})},
we get
\begin{eqnarray*}
{\mathbf{E}}[\|X^{\varepsilon }_{i,s}\|_\beta^p]&\leq&
{\mathbf{E}}[Z_i^p]\leq2^{p-1}
({\mathbf{E}}[|Z_i-\|\mathcaligr{S}_sx_i\|_\beta
|^p]+\|\mathcaligr{S}
_sx_i\|_\beta^p)
\\
&\leq& C_p\bigl((\|\mathcaligr{S}_sx_1\|_\beta \|\mathcaligr{S}_sx_2\|
_\beta)^{p/2}
+\|\mathcaligr{S}_sx_i\|_\beta^p\bigr)
\\
&\leq& C_p e^{pMs}\bigl((\|x_1\|_\beta \|x_2\|_\beta
)^{p/2}+\|
x_i\|_\beta^p\bigr).
\end{eqnarray*}
This shows \mbox{\rm(\ref{E5.2})} and completes the proof of
Theorem~\ref{T1}.

\subsection[Proof of Theorem 2]{Proof of Theorem~\protect\ref{T2}}
\label{S5.2}
Theorem~\ref{T2} is a direct consequence of Theorem~\ref{T1},
Corollary~\ref{C3.2.3} and \mbox{\rm(\ref{E2.1.2})}.

\section*{Acknowledgment}
The authors wish to thank an anonymous referee for valuable comments.

%

\printaddresses


\begin{thebibliography}{13}

\bibitem{Aldous1978}
\begin{barticle}[msn]
\bauthor{\bsnm{Aldous},~\bfnm{David}\binits{D.}}
(\byear{1978}).
\btitle{Stopping times and tightness}.
\bjournal{Ann. Probability}
\bvolume{6}
\bpages{335--340}.
\bmrnumber{MR0474446}
\end{barticle}
\endbibitem

\bibitem{CoxKlenke2000}
\begin{barticle}[msn]
\bauthor{\bsnm{Cox},~\bfnm{J.~Theodore}\binits{J.~T.}} \AND
  \bauthor{\bsnm{Klenke},~\bfnm{Achim}\binits{A.}}
(\byear{2000}).
\btitle{Recurrence and ergodicity of interacting particle systems}.
\bjournal{Probab. Theory Related Fields}
\bvolume{116}
\bpages{239--255}.
\bmrnumber{MR1743771}
\end{barticle}
\endbibitem

\bibitem{CoxKlenkePerkins2000}
\begin{bincollection}[msn]
\bauthor{\bsnm{Cox},~\bfnm{J.~Theodore}\binits{J.~T.}},
  \bauthor{\bsnm{Klenke},~\bfnm{Achim}\binits{A.}} \AND
  \bauthor{\bsnm{Perkins},~\bfnm{Edwin~A.}\binits{E.~A.}}
(\byear{2000}).
\btitle{Convergence to equilibrium and linear systems duality}.
In \bbooktitle{Stochastic Models}
(\beditor{L. G. Gorostiza and B. G. Ivanoff}, eds.)
\bpages{41--66}.
\bpublisher{Amer. Math. Soc.}, \baddress{Providence, RI}.
\bmrnumber{MR1765002}
\end{bincollection}
\endbibitem

\bibitem{DawsonPerkins1998}
\begin{barticle}[msn]
\bauthor{\bsnm{Dawson},~\bfnm{Donald~A.}\binits{D.~A.}} \AND
  \bauthor{\bsnm{Perkins},~\bfnm{Edwin~A.}\binits{E.~A.}}
(\byear{1998}).
\btitle{Long-time behavior and coexistence in a mutually catalytic branching
  model}.
\bjournal{Ann. Probab.}
\bvolume{26}
\bpages{1088--1138}.
\bmrnumber{MR1634416}
\end{barticle}
\endbibitem

\bibitem{DellacherieMeyer1983}
\begin{bbook}[msn]
\bauthor{\bsnm{Dellacherie},~\bfnm{Claude}\binits{C.}} \AND
  \bauthor{\bsnm{Meyer},~\bfnm{Paul-Andr{\'e}}\binits{P.-A.}}
(\byear{1980}).
\btitle{Probabilit\'es et Potentiel. {C}hapitres {V} \`a {VIII}},
\bedition{Revised} ed.
\bseries{Actualit\'es Scientifiques et Industrielles [Current Scientific and
  Industrial Topics]}
\bvolume{1385}.
\bpublisher{Hermann}, \baddress{Paris}.
\bmrnumber{MR566768}
\end{bbook}
\endbibitem

\bibitem{EtheridgeFleischmann2004}
\begin{barticle}[msn]
\bauthor{\bsnm{Etheridge},~\bfnm{Alison~M.}\binits{A.~M.}} \AND
  \bauthor{\bsnm{Fleischmann},~\bfnm{Klaus}\binits{K.}}
(\byear{2004}).
\btitle{Compact interface property for symbiotic branching}.
\bjournal{Stochastic Process. Appl.}
\bvolume{114}
\bpages{127--160}.
\bmrnumber{MR2094150}
\end{barticle}
\endbibitem

\bibitem{EthierKurtz1986}
\begin{bbook}[msn]
\bauthor{\bsnm{Ethier},~\bfnm{Stewart~N.}\binits{S.~N.}} \AND
  \bauthor{\bsnm{Kurtz},~\bfnm{Thomas~G.}\binits{T.~G.}}
(\byear{1986}).
\btitle{Markov Processes}.
\bpublisher{Wiley}, \baddress{New York}.
\bmrnumber{MR838085}
\end{bbook}
\endbibitem

\bibitem{JacodShiryaev1987}
\begin{bbook}[msn]
\bauthor{\bsnm{Jacod},~\bfnm{Jean}\binits{J.}} \AND
  \bauthor{\bsnm{Shiryaev},~\bfnm{Albert~N.}\binits{A.~N.}}
(\byear{1987}).
\btitle{Limit Theorems for Stochastic Processes}.
\bseries{Grundlehren der Mathematischen Wissenschaften [Fundamental Principles
  of Mathematical Sciences]}
\bvolume{288}.
\bpublisher{Springer}, \baddress{Berlin}.
\bmrnumber{MR959133}
\end{bbook}
\endbibitem

\bibitem{KM1}
\begin{bmisc}[unstr]
\bauthor{\bsnm{Klenke},~\bfnm{A.}\binits{A.}} \AND
  \bauthor{\bsnm{Mytnik},~\bfnm{L.}\binits{L.}}
 (2008). Infinite rate mutually catalytic
  branching. Preprint. Available at arXiv:\href{http://arxiv.org/abs/0809.4554}{0809.4554} [math.PR].
\end{bmisc}
\endbibitem

\bibitem{KM2}
\begin{bmisc}[unstr]
\bauthor{\bsnm{Klenke},~\bfnm{A.}\binits{A.}} \AND
  \bauthor{\bsnm{Mytnik},~\bfnm{L.}\binits{L.}}
  (2008). Infinite rate mutually catalytic
  branching in infinitely many colonies. {C}onstruction, characterization and
  convergence. Preprint. Available at arXiv:\href{http://arxiv.org/abs/0901.0623}{0901.0623}
  [math.PR].
\end{bmisc}
\endbibitem

\bibitem{KM3}
\begin{bmisc}[unstr]
\bauthor{\bsnm{Klenke},~\bfnm{A.}\binits{A.}} \AND
  \bauthor{\bsnm{Mytnik},~\bfnm{L.}\binits{L.}}
   (2009). Infinite rate mutually catalytic
  branching in infinitely many colonies. The longtime behaviour. Preprint.
\end{bmisc}
\endbibitem

\bibitem{Liggett1985}
\begin{bbook}[msn]
\bauthor{\bsnm{Liggett},~\bfnm{Thomas~M.}\binits{T.~M.}}
(\byear{1985}).
\btitle{Interacting Particle Systems}.
\bseries{Grundlehren der Mathematischen Wissenschaften [Fundamental Principles
  of Mathematical Sciences]}
\bvolume{276}.
\bpublisher{Springer}, \baddress{New York}.
\bmrnumber{MR776231}
\end{bbook}
\endbibitem

\bibitem{Oeler2008}
\begin{bmisc}[unstr]
\bauthor{\bsnm{Oeler},~\bfnm{Mario}\binits{M.}}
(\byear{2008}).
 {Mutually catalytic branching at infinite
  rate}. Ph.D. thesis, Univ. Mainz.
\end{bmisc}
\endbibitem

\end{thebibliography}
\end{document}